\documentclass[11pt,a4paper]{amsart}
\usepackage{amssymb,amsmath}
\usepackage{mathpazo}

\headsep=1cm \numberwithin{equation}{section}
\newtheorem{theorem}{Theorem}[section]
\newtheorem{lemma}[theorem]{Lemma}
\newtheorem{proposition}[theorem]{Proposition}
\newtheorem{corollary}[theorem]{Corollary}

\theoremstyle{definition}

\theoremstyle{remark}
\newtheorem{remark}[theorem]{Remark}
\newtheorem{example}[theorem]{Example}

\newcommand{\Ass}{\operatorname{Ass}}

\newcommand{\Assh}{\operatorname{Assh}}

\newcommand{\cd}{\operatorname{cd}}

\newcommand{\Ht}{\operatorname{ht}}

\newcommand{\Var}{\operatorname{Var}}

\newcommand{\Supp}{\operatorname{Supp}}

\newcommand{\Att}{\operatorname{Att}}
\newcommand{\Ann}{\operatorname{Ann}}
\newcommand{\Rad}{\operatorname{Rad}}

\newcommand{\Max}{\operatorname{Max}}

\newcommand{\lo}{\longrightarrow}
\newcommand{\fm}{\frak{m}}
\newcommand{\fp}{\frak{p}}

\newcommand{\fa}{\frak{a}}
\newcommand{\fb}{\frak{b}}
\newcommand{\fc}{\frak{c}}

\newcommand{\fn}{\frak{n}}
\newcommand{\fM}{\frak{M}}
\newcommand{\fX}{\frak{X}}
\newcommand{\fW}{\frak{W}}

\newcommand{\fZ}{\frak{Z}}

\begin{document}
\author[K. Divaani-Aazar]{Kamran Divaani-Aazar}
\title[Vanishing of the top local cohomology ...]
{Vanishing of the top local cohomology modules over  Noetherian
rings}

\address{K. Divaani-Aazar, Department of Mathematics, Az-Zahra University,
Vanak, Post Code 19834, Tehran, Iran.} \email{kdivaani@ipm.ir}



\begin{abstract}Let $R$ be a (not necessarily local) Noetherian ring
and $M$ a finitely generated $R$-module of finite dimension $d$.
Let $\fa$ be an ideal of $R$ and $\fM$ denote the intersection of
all prime ideals $\fp\in \Supp_RH^d_{\fa}(M)$. It is shown that
$$H^d_{\fa}(M)\simeq H^d_{\fM}(M)/\displaystyle{\sum_{n\in
\mathbb{N}}}<\fM>(0:_{H^d_{\fM}(M)}\fa^n),$$ where for an Artinian
$R$-module $A$ we put $<\fM>A=\cap_{n\in \mathbb{N}} \fM^nA$.  As
a consequence, it is proved that for all ideals $\fa$ of $R$,
there are only finitely many non-isomorphic top local cohomology
modules $H^d_{\fa}(M)$ having the same support. In addition, we
establish an analogue of the Lichtenbaum-Hartshorne Vanishing
Theorem over rings that need not be local.
\end{abstract}

\maketitle

\section{Introduction}

Throughout this paper, let R denote a commutative Noetherian ring.
Let $M$ be a finitely generated $R$-module of finite dimension $d$
and $\fa$ an ideal of $R$. The present article is concerned with
the top local cohomology module $H^d_{\fa}(M)$. We refer the
reader to [{\bf  3}] for more details about local cohomology. By
Grothendieck's Vanishing Theorem [{\bf 3}, Theorem 6.1.2], it is
known that $H^i_{\fa}(M)=0$ for all $i>\dim M$. So $H^d_{\fa}(M)$
is the last possible non-vanishing local cohomology module of $M$.
Also, by [{\bf  3}, Exercise 7.1.7] the top local cohomology
module $H^d_{\fa}(M)$ is Artinian. There are many papers
concerning the top local cohomology modules of finitely generated
modules over local rings. But, according to the best knowledge of
the author, [{\bf 2}] and [{\bf 4}] are the only existing articles
studying such local cohomology modules over general Noetherian
rings. In this paper, we investigate the structure of the top
local cohomology modules of finitely generated modules over rings
that need not be local.

When $R$ is local with the maximal ideal $\fm$, it is proved that
there is a natural isomorphism $H^d_{\fa}(M)\simeq
H^d_{\fm}(M)/\Sigma_{n\in
\mathbb{N}}<\fm>(0:_{H^d_{\fm}(M)}\fa^n)$, see [{\bf 10}, Theorem
3.2]. As a result, in [{\bf  10}] a new proof is provided for the
Lichtenbaum-Hartshorne Vanishing Theorem. In Section 2, we
establish an analogue of the above isomorphism over rings that are
not necessarily local. To be more precise, we will prove that if
$\fM$ denotes the intersection of all prime ideals $\fp\in
\Supp_RH^d_{\fa}(M)$, then there is a natural isomorphism
$$H^d_{\fa}(M)\simeq H^d_{\fM}(M)/\displaystyle{\sum_{n\in
\mathbb{N}}}<\fM>(0:_{H^d_{\fM}(M)}\fa^n).$$ This will be proved
in Theorem 2.3.

Knowing more about $\Att_RH^d_{\fa}(M)$, the set of attached primes of
$H^d_{\fa}(M)$, could lead to better understanding of the
structure of the top local cohomology module $H^d_{\fa}(M)$. In
particular, knowing $\Att_RH^d_{\fa}(M)$ implies vanishing results
for $H^d_{\fa}(M)$. In the case $R$ is local, the set
$\Att_RH^d_{\fa}(M)$ is already determined (see e.g. [{\bf 18}],
[{\bf 10}] and [{\bf  6}]). In Theorem 2.5 below, we determine the
set $\Att_RH^d_{\fa}(M)$ without the assumption that $R$ is local,
namely we show that
$$\Att_RH^d_{\fa}(M)=\{\fp\in \Assh_RM:
\cd_R(\fa,R/\fp)=d\}$$(here for an $R$-module $N$, $\cd_R(\fa,N)$
denotes the cohomological dimension of $N$ with respect to the
ideal $\fa$). Then as an application, we provide an improvement of the main result
of [{\bf 2}].  Next, for a finitely generated $R$-module $N$ so that $H^c_{\fa}(N)$, $c:=\cd_R(\fa,N)$, is representable, we examine the set $\Att_RH^c_{\fa}(N)$.

In Section 3, first we show that for all ideals $\fa$ of $R$,
there are only finitely many non-isomorphic top local cohomology
modules $H^d_{\fa}(M)$ having the same support. Next, as an
application of Theorems 2.3 and 2.5, we extend the
Licthenbum-Hartshorne Vanishing Theorem to (not necessarily local)
Noetherian rings. Namely, we prove that if $\fM$ is as above and
$T$ denotes the $\fM$-adic completion of $R$, then
the following are equivalent:\\
i) $H^d_{\fa}(M)=0$.\\
ii) $H^d_{\fM}(M)=\displaystyle{\sum_{n\in
\mathbb{N}}}<\fM>(0:_{H^d_{\fM}(M)}\fa^n)$.\\
iii) For any integer $l\in \mathbb{N}$, there exists an $n=n(l)\in
\mathbb{N}$ such that $$0:_{H^d_{\fM}(M)}\fa^l\subseteq
<\fM>(0:_{H^d_{\fM}(M)}\fa^n).$$\\
iv) $\dim T/\fa T+\fp>0$ for all
$\fp\in \Assh_T(M\otimes_RT)$.\\
v) $\cd_R(\fa,R/\fp)<d$ for all $\fp\in \Assh_RM$.

Throughout the paper, for an $R$-module $M$, $\Assh_RM$ denotes
the set of all associated prime ideals $\fp$ of $M$ such that
$\dim R/\fp=\dim M$. Also, for an Artinian $R$-module $A$, we
denote $\cap_{n\in \mathbb{N}}\fa^nA$ by $<\fa>A$.

\section{Attached prime ideals}

A nonzero $R$-module $S$ is called {\it secondary} if for each
$x\in R$ the multiplication map induced by $x$ on $S$ is either
surjective or nilpotent. If $S$ is secondary, then the ideal
$\fp:=\Rad(\Ann_RS)$ is a prime ideal and $S$ is called
$\fp$-secondary. For an $R$-module $M$, a secondary representation
of $M$ is an expression for $M$ as a sum of finitely many
secondary submodules of $M$. An $R$-module $M$ is said to be {\it
representable} if it has a secondary representation. From any
secondary representation for an $R$-module $M$, one can obtain
another one as $M=S_1+\dots +S_n$ such that the prime ideals
$\fp_i:=\Rad(\Ann_RS_i), i=1,\dots , n$ are all distinct and
$S_j\nsubseteq \Sigma_{i\neq j}S_i$ for all $j=1,\dots ,n$. A such
secondary representation for $M$ is said to be minimal. It is
shown that the set $\{\fp_1,\dots , \fp_n\}$ is independent of the
chosen minimal secondary representation for $M$. This set is
denoted by $\Att_RM$ and each element of this set is said to be an
attached prime ideal of $M$. It is known that a representable
$R$-module $M$ is zero if and only if $\Att_RM=\emptyset$ and that
if $0\lo N\lo M\lo L\lo 0$ is an exact sequence of representable
$R$-modules and $R$-homomorphisms, then $\Att_RL\subseteq
\Att_RM\subseteq \Att_RN\cup \Att_RL$. Also, it is known that any
Artinian $R$-module is representable. For more information about
the theory of secondary representations see [{\bf  12}] or [{\bf
14}, Section 6, Appendix].

\begin{lemma} i) Let $f:R\lo U$ be a ring homomorphism
and $M$ a representable $U$-module. Then $M$ is also representable
as an $R$-module and
$\Att_RM=\{f^{-1}(\fp):\fp\in \Att_UM\}$.\\
ii) Let $A$ be an Artinian $R$-module. Then $\Supp_RA$ equals
$\Ass_RA$ and is a finite subset of $\Max R$. Moreover, if
$\Supp_RA=\{\fm_1,\dots ,\fm_t\}$, then the natural
$R$-homomorphism $\psi:A\lo \oplus_{i=1}^tA_{\fm_i}$ is an
isomorphism. In particular, $\Att_RA=\displaystyle{\bigcup_{i=1}^t}
\Att_RA_{\fm_i}$.\\
iii) Let $\fm_1,\dots ,\fm_t$ be distinct maximal ideals of $R$
and $A_1,\dots ,A_t$ Artinian $R$-modules so that
$\Supp_RA_i=\{\fm_i\}$ for all $i=1,\dots ,t$. Let
$A=\oplus_{i=1}^tA_i$. Then for any ideal $\fa$ of $R$ such that
$\fa\subseteq \fM:=\cap_{i=1}^t\fm_i$, there is a natural
isomorphism
$$\dfrac{A}{\displaystyle{\sum_{n\in
\mathbb{N}}}<\fM>(0:_A\fa^n)}\simeq
\displaystyle{\bigoplus_{i=1}^t}\dfrac{A_i}{\displaystyle{\sum_{n\in
\mathbb{N}}}<\fm_i>(0:_{A_i}\fa^n)}.$$
\end{lemma}

{\bf Proof.} i) holds by [{\bf 15}, Proposition 4.1].

ii) The first assertion of (ii) holds by [{\bf  17}, Exercises
8.49 and 9.43]. Now, we are going to prove the second assertion of
(ii). It follows by [{\bf  17}, Exercise 8.49], that
$A=\displaystyle{\oplus_{i=1}^t}\Gamma_{\fm_i}(A)$. This yields
that for each i, $A_{\fm_i}\simeq \Gamma_{\fm_i}(A)$, and so
$A_{\fm_i}$, as an $R$-module, supported only at the maximal ideal
$\fm_i$. So $\psi_{\fm}:A_{\fm}\lo
(\displaystyle{\oplus_{i=1}^t}A_{\fm_i})_{\fm}$ is an isomorphism
for any maximal ideal $\fm$ of $R$. Thus $\psi$ is an isomorphism,
as claimed. Finally, the last assertion of (ii) is immediate by
(i) and the fact that for any given finitely many secondary
representable $R$-modules $M_1,\dots ,M_t$, it turns out that
$\displaystyle{\oplus_{i=1}^t}M_i$ is also representable and that
$$\Att_R(\displaystyle{\bigoplus_{i=1}^t}M_i)=\bigcup_{i=1}^t\Att_RM_i.$$

iii) First of all note that for any Artinian $R$-module $B$ and
any two ideals $\fa$, $\fb$ of $R$, it is easy to see that
$\{\dfrac{0:_B\fa^n}{<\fb>(0:_B\fa^n)}\}_{n\in \mathbb{N}}$, with
the natural maps induced by the identity map of $B$, is a direct
system and that $$\dfrac{\displaystyle{\sum_{n\in
\mathbb{N}}}(0:_B\fa^n)}{\displaystyle{\sum_{n\in
\mathbb{N}}}<\fb>(0:_B\fa^n)}$$ is its direct limit. In
particular, if $\fa\subseteq \cap_{\fm\in \Supp_RB}\fm$, then each
element of $B$ is annihilated by some power of $\fa$, and so
$$\underset{n}{\varinjlim}\frac{0:_B\fa^n}{<\fb>(0:_B\fa^n)}=\frac{B}
{\displaystyle{\sum_{n\in
\mathbb{N}}}<\fb>(0:_B\fa^n)}.$$ Next, note that $A_{\fm_i}\simeq
A_i$ for all $i=1,\dots ,t$. Thus in view of (ii), we have the
following isomorphisms

$$\begin{array}{llll} \dfrac{A}{\displaystyle{\sum_{n\in
\mathbb{N}}}<\fM>(0:_A\fa^n)}
& \simeq
\underset{n}{\varinjlim}\dfrac{0:_A\fa^n}{<\fM>(0:_A\fa^n)}\\
& \simeq \underset{n}{\varinjlim}[(\dfrac{0:_A\fa^n}{<\fM>
(0:_A\fa^n)})_{\fm_1}\oplus
\dots \oplus(\dfrac{0:_A\fa^n}{<\fM>(0:_A\fa^n)})_{\fm_t}]\\
& \simeq
\underset{n}{\varinjlim}[\dfrac{0:_{A_1}\fa^n}{<\fm_1>(0:_{A_1}\fa^n)}
\oplus \dots \oplus \dfrac{0:_{A_t}\fa^n}{<\fm_t>(0:_{A_t}\fa^n)}]\\
& \simeq
\displaystyle{\bigoplus_{i=1}^t}[\underset{n}{\varinjlim}\dfrac{0:_{A_i}
\fa^n}{<\fm_i>(0:_{A_i}\fa^n)}]\\
& \simeq
\displaystyle{\bigoplus_{i=1}^t}\dfrac{A_i}{\displaystyle{\sum_{n\in
\mathbb{N}}} <\fm_i>(0:_{A_i}\fa^n)}. \Box \end{array}$$

\begin{remark} i) Let $\fa$ be an ideal of $R$. For a prime
ideal $\fp$ of $R$, we say that $\fa$ is {\it formally isolated
at} $\fp$ if $\fa\subseteq \fp$ and if there is some prime ideal
$\fp^*$ of $\hat{R_{\fp}}$ such that $\dim
\hat{R_{\fp}}/\fp^*=\Ht(\fp)$ and that $\dim
\hat{R_{\fp}}/\fa\hat{R_{\fp}}+\fp^*=0$. Assume that $R$ has
finite dimension $d$, and let $\mathcal{P}_{\fa}$ denote the set
of all prime ideals $\fp$ such that $\Ht(\fp)=d$ and such that
$\fa$ is formally isolated at $\fp$. Then, by [{\bf 2}, Theorem
3.3 (b)] for any finitely generated faithful $R$-module $M$, we
have $\Supp_RH^d_{\fa}(M)=\mathcal{P}_{\fa}$.\\
ii) Let $M$ be a finitely generated $R$-module of finite dimension
$d$. Let $\mathcal{P}_{\fa,M}$ denote the set of all $\fp\in \Var(
\Ann_RM+\fa)$ so that there is some prime $\fp^*\in
\Supp_{\hat{R_{\fp}}}\hat{M_{\fp}}$ such that $\dim
\hat{R_{\fp}}/\fp^*=d$ and that $\dim \hat{R_{\fp}}/\fa
\hat{R_{\fp}}+\fp^*=0$. Then, by adapting the method of the proof
of [{\bf 2}, Theorem 3.3(b)], one can easily deduce that
$\Supp_RH^d_{\fa}(M)=\mathcal{P}_{\fa,M}$. Also, in Corollary 4.1
below,  we establish another characterization of
$\mathcal{P}_{\fa,M}$.
\end{remark}

In the remainder of the paper, for a finitely generated $R$-module
$M$ of finite dimension $d$ and an ideal $\fa$ of $R$, let
$\mathcal{P}_{\fa,M}$ be as in Remark 2.2 (ii).

\begin{theorem} Let $\fa$ be an ideal of $R$, $M$ a finitely
 generated $R$-module of finite dimension $d$ and
 $\fM=\displaystyle{\bigcap_{\fp\in \mathcal{P}_{\fa,M}}}\fp$. There is a
 natural isomorphism $$H^d_{\fa}(M)\simeq H^d_{\fM}(M)/
\displaystyle{\sum_{n\in
\mathbb{N}}}<\fM>(0:_{H^d_{\fM}(M)}\fa^n).$$
\end{theorem}

{\bf Proof.}  By Remark 2.2 (ii), we have
$\Supp_RH^d_{\fa}(M)=\mathcal{P}_{\fa,M}$. Let
$\Supp_RH^d_{\fa}(M)=\{\fm_1,\dots ,\fm_t\}$ and for each $i$
denote the local ring $R_{\fm_i}$ by $R_i$.

Let $\fa$ be an ideal of a local ring $(U,\fn)$. By [{\bf  10},
Theorem 3.2], it turns out that for any finitely generated
$U$-module $M$, there is a natural isomorphism
$$H^d_{\fa}(M)\simeq H^d_{\fn}(M)/\displaystyle{\sum_{n\in
\mathbb{N}}}<\fn>(0:_{H^d_{\fn}(M)}\fa^n),$$ where $d=\dim M$.
Observe that by the Flat Base Change Theorem [{\bf  3}, Theorem
4.3.2] and Lemma 2.1 (ii) the modules $H^d_{\fm_iR_i}(M_{\fm_i})$
and $H^d_{\fm_i}(M)$ are isomorphic for all $1\leq i\leq t$.
Therefore applying Lemma 2.1 (ii) again, provides the following
isomorphisms
$$\begin{array}{llll} H^d_{\fa}(M)& \simeq
\displaystyle{\bigoplus_{i=1}^t}H^d_{\fa R_i}(M_{\fm_i})\\
& \simeq
\displaystyle{\bigoplus_{i=1}^t}\dfrac{H^d_{\fm_iR_i}(M_{\fm_i})}{
\displaystyle{\sum_{n\in \mathbb{N}}}<\fm_iR_i>(0:_{H^d_{\fm_iR_i}
(M_{\fm_i})}\fa^nR_i)}\\
& \simeq
\displaystyle{\bigoplus_{i=1}^t}\dfrac{H^d_{\fm_i}(M)}{\displaystyle
{\sum_{n\in
\mathbb{N}}}<\fm_i>(0:_{H^d_{\fm_i}(M)}\fa^n)}.\end{array}$$ On
the other hand, the Mayer-Vietoris sequence for local cohomology
[{\bf 3}, Theorem 3.2.3] yields the following isomorphism
$$H^d_{\fM}(M)\simeq \displaystyle{\bigoplus_{i=1}^t}H^d_{\fm_i}(M).$$
This finishes the proof, by Lemma 2.1 (iii). $\Box$

Recall that for an $R$-module $M$, the cohomological dimension of
$M$ with respect to an ideal $\fa$ of $R$ is defined as
$\cd_R(\fa,M):=\sup \{i\in \mathbb{N}_0: H^i_{\fa}(M)\neq 0\}$. It
is appropriate to list some basic properties of this notion. First
of all note that, it is immediate by Grothendieck's Vanishing
Theorem, that $\cd_R(\fa,M)\leq \dim M$. Next, note that if $V$ is
a multiplicative subset of $R$, then it becomes clear by the Flat
Base Change Theorem, that $\cd_{V^{-1}R}(\fa V^{-1}R, V^{-1}M)\leq
\cd_R(\fa,M)$. Also, if $M$ and $L$ are two finitely generated
$R$-modules so that $\Supp_RL\subseteq \Supp_RM$, then [{\bf  9},
Theorem 2.2] implies that $\cd_R(\fa,L)\leq \cd_R(\fa,M)$. For
further details concerning this notion, we refer the reader to
[{\bf  11}] and [{\bf  9}].

\begin{lemma} Let $\fa$ be an ideal of a local ring
$(R,\fm)$ and $d$ a natural number. For any prime ideal $\fp$ of
$R$ so that $\dim R/\fp\leq d$,
the following are equivalent:\\
i) $\cd_R(\fa,R/\fp)=d$.\\
ii) $\fp$ is the contraction to $R$ of a prime ideal $\fp^*$ of
$\hat{R}$ such that $\dim \hat{R}/\fp^*=d$ and $\dim
\hat{R}/\fa\hat{R}+\fp^*=0$.
\end{lemma}

{\bf Proof.} Let $M$ be a finitely generated $R$-module of
dimension $d$. Then by the Lichtenbaum-Hartshorne Vanishing
Theorem, it turns out that $H^d_{\fa}(M)\neq 0$ if and only if
there exists $\fp^*\in \Assh_{\hat{R}}\hat{M}$ such that $\dim
\hat{R}/\fa\hat{R}+\fp^*=0$ (see e.g. [{\bf  10}, Corollary 3.4]).
Assume that (i) holds. Then $H^d_{\fa}(R/\fp)\neq 0$, and so there
exists $\fp^*\in \Assh_{\hat{R}}(\hat{R}/\fp\hat{R})$ such that
$\dim \hat{R}/\fa\hat{R}+\fp^*=0$. Since $H^d_{\fa}(R/\fp)\neq 0$,
by Grothendieck's Vanishing Theorem, we have $\dim R/\fp=d$. Thus
$$\dim \hat{R}/\fp^*=\dim \hat{R}/\fp\hat{R}=d.$$ On the other hand,
by [{\bf 14}, Theorem 23.2 (i)], we have
$$\{\fp\}=\Ass_R(R/\fp)=\{Q\cap R: Q\in
\Ass_{\hat{R}}(\hat{R}/\fp\hat{R})\}.$$ Hence $\fp=\fp^*\cap R$,
and so (ii) follows.

Now, assume that (ii) holds. We have $$d\geq \dim R/\fp= \dim
\hat{R}/\fp\hat{R}\geq \dim \hat{R}/\fp^*=d.$$ So $\dim R/\fp=d$.
In particular, $\fp^*$ is minimal over $\fp \hat{R}$, and so
$\fp^*\in \Assh_{\hat{R}}(\hat{R}/\fp\hat{R})$. Thus
$H^d_{\fa}(R/\fp)\neq 0$, by the Lichtenbaum-Hartshorne Vanishing
Theorem (that we commented earlier its statement in the beginning
of the proof). Therefore $\cd_R(\fa,R/\fp)=d$, as required. $\Box$

The following extends the main result of [{\bf  6}] to general
Noetherian rings.

\begin{theorem} (See [{\bf 4}, Theorem 1.2]) Let $\fa$ be an
ideal of $R$ and $M$ a finitely generated $R$-module of finite
dimension $d$. Then
$$\Att_RH^d_{\fa}(M)=\{\fp\in \Assh_RM: \cd_R(\fa,R/\fp)=d\}.$$
\end{theorem}

{\bf Proof.} Assume that $\Supp_RH^d_{\fa}(M)=\{\fm_1,\dots
,\fm_t\}$. Then by the Flat Base Change Theorem and Lemma 2.1
(ii), it follows that
$$\Att_RH^d_{\fa}(M)=\bigcup_{i=1}^t\Att_RH^d_{\fa
R_{\fm_i}}(M_{\fm_i}).$$ In the remainder of the proof, we will
use this equality without further comment.

Let $M$ be a finitely generated module over a local ring
$(U,\fn)$. Then by [{\bf  10}, Corollary 3.3] for any ideal $\fa$
of $U$, $\Att_{\hat{U}}H^{\dim M}_{\fa}(M)$ consists of all
$\fp\in \Assh_{\hat{U}}\hat{M}$ such that $\dim
\hat{U}/\fa\hat{U}+\fp=0$. Fix $1\leq i\leq t$. Since
$$H^d_{\fa R_{\fm_i}}(M_{\fm_i})\simeq (H^d_{\fa}(M))_{\fm_i}\neq
0,$$ we have $\dim M_{\fm_i}=d$. It now follows, by Lemma 2.1 (i)
and Lemma 2.4 that

$$\begin{array}{llll} \Att_{R_{\fm_i}}H^d_{\fa R_{\fm_i}}(M_{\fm_i})&
=\{Q\cap R_{\fm_i}:Q\in \Assh_{\hat{R}_{\fm_i}}\hat{M}_{\fm_i} ,
\dim
\hat{R}_{\fm_i}/\fa\hat{R}_{\fm_i}+Q=0\} \\
& =\{\fp R_{\fm_i}\in
\Assh_{R_{\fm_i}}M_{\fm_i}:\cd_{R_{\fm_i}}(\fa
R_{\fm_i},R_{\fm_i}/ \fp R_{\fm_i})=d\}.\end{array}$$
Because
$\dim M_{\fm_i}=\dim M=d$ and
$$\Ass_{R_{\fm_i}}M_{\fm_i}=\{\fp R_{\fm_i}:\fp\subseteq \fm_i
\ and \  \fp\in \Ass_RM\},$$ it follows that
$\Assh_{R_{\fm_i}}M_{\fm_i}$ consists of all prime ideals $\fp
R_{\fm_i}\in \Ass_{R_{\fm_i}}M_{\fm_i}$ such that $\fp\in
\Assh_RM$. Hence, if $\fp\in \Att_RH^d_{\fa}(M)$, then $\fp\in
\Assh_RM$ and $\cd_R(\fa,R/\fp)=d$.

Conversely, assume that $\fp\in \Assh_RM$ is such that
$\cd_R(\fa,R/\fp)=d$. Let $\fm\in \Supp_RH^d_{\fa}(R/\fp)$. Then
$H^d_{\fa R_{\fm}}(R_{\fm}/\fp R_{\fm})\neq 0$, and so $\dim
R_{\fm}/\fp R_{\fm}=d$. Hence, we have $\cd_{R_{\fm}}(\fa
R_{\fm},R_{\fm}/\fp R_{\fm})=d$ and $\fp R_{\fm}\in
\Assh_{R_{\fm}}M_{\fm}$. By Lemma 2.4, $\fp R_{\fm}$ is the
contraction to $R_{\fm}$ of a prime ideal $\fp^*$ of
$\hat{R}_{\fm}$ such that $\dim \hat{R}_{\fm}/\fp^*=d$ and $\dim
\hat{R}_{\fm}/\fa \hat{R}_{\fm}+\fp^*=0$. It is easy to see that
$\fp^*\in \Assh_{\hat{R}_{\fm}}\hat{M}_{\fm}$, and so by Lemma 2.1
(i) and the above mentioned result of [{\bf 10}], it turns out
that $\fp R_{\fm}\in \Att_{R_{\fm}}H^d_{\fa R_{\fm}}(M_{\fm})$.
Hence $\fp\in \Att_RH^d_{\fa}(M)$, by using Lemma 2.1 (i) again.
Note that, since $\Att_RH^d_{\fa}(M)_{\fm}$ is not empty, it
follows that $\fm\in \Supp_RH^d_{\fa}(M)$. $\Box$

\begin{example} In [{\bf  8}, Corollary 3.3], the fact that the
top local cohomology modules of finitely generated modules of
finite dimension are Artinian is extended to an strictly larger
class of modules. Namely, it is shown that if $\fa$ is an ideal of
$R$ and $M$ a $ZD$-module of finite dimension $d$ such that
$\fa$-relative Goldie dimension of any quotient of $M$ is finite,
then $H^d_{\fa}(M)$ is Artinian. It would be interesting to know
whether the conclusion of Theorem 2.5 remains valid for this
larger class of modules. Unfortunately, this is not the case, even
if $R$ is local. To this end, let $(R,\fm)$ be a local ring with
$\dim R>0$. Take $\fa=\fm$ and $M=E(R/\fm)$, the injective envelop
of the residue field of $R$. Then $M$ is a $ZD$-module and
$\fa$-relative Goldie dimension of any quotient of $M$ is finite.
We have $$\Att_RH^0_{\fa}(M)=\Att_RM=\Ass_RR,$$ while the maximal
ideal $\fm$ is the only element of the set $$\{\fp\in \Assh_RM:
\cd_R(\fa,R/\fp)=0\}.$$
\end{example}

As a corollary to Theorem 2.5, we present an improvement of
the main result of [{\bf  2}]. In the sequel, let
$\mathcal{P}_{\fa}$ be as in Remark 2.2 (i).

\begin{corollary} Let $\fa$ and $\fb$ be two ideals of $R$ and assume
that $R/\fb$ has finite dimension $d$. Then\\
i) $\mathcal{P}_{\fa,R/\fb}=\{\fm\in \Max R:\exists \fp\in
\Assh_R(R/\fb) \ such \ that \  \fp\subseteq \fm \ and \
\cd_{R_{\fm}}(\fa R_{\fm},R_{\fm}/\fp R_{\fm})\\=d\}$. In
particular, if $R$ has finite dimension $d$, then
$$\mathcal{P}_{\fa}=\{\fm\in \Max R:\exists \fp\in \Assh_RR
\ such \  that \  \fp\subseteq \fm \ and \ \cd_{R_{\fm}}
(\fa R_{\fm},R_{\fm}/\fp R_{\fm})=d\}.$$\\
ii) For any finitely generated $R$-module $M$ so that
$\Assh_RM=\Assh_R(R/\fb)$, we have
$\Supp_RH^d_{\fa}(M)=\mathcal{P}_{\fa,R/\fb}$. In particular,
$\mathcal{P}_{\fa,R/\fb}$ is a finite set.\\
iii) (See [{\bf 4}, Theorem 1.3 (g)]) If $d>0$, then for any $M$
as in (ii), the $R_{\fm}$-module $(H^d_{\fa}(M))_{\fm}$ is not
finitely generated for all $\fm\in \mathcal{P}_{\fa,R/\fb}$.
\end{corollary}

{\bf Proof.} First, it should be noted that
$\mathcal{P}_{\fa}=\mathcal{P}_{\fa,R}$. By Remark 2.2 (ii), we
have $\Supp_RH^d_{\fa}(R/\fb)=\mathcal{P}_{\fa,R/\fb}$. Hence, to
prove (i) and (ii), it will be enough to show that for any
finitely generated $R$-module $M$ with $\Assh_RM=\Assh_R(R/\fb)$,
$\Supp_RH^d_{\fa}(M)$ consists of all maximal ideals $\fm$ of $R$
so that there exists a prime ideal $\fp\in \Assh_R(R/\fb)$ such
that $\fp\subseteq \fm $ and $\cd_{R_{\fm}}(\fa
R_{\fm},R_{\fm}/\fp R_{\fm})=d$. Assume that $M$ is a finitely
generated $R$-module with $\Assh_RM=\Assh_R(R/\fb)$, and let
$\fm\in \Supp_RH^d_{\fa}(M)$. Then $H^d_{\fa R_{\fm}}(M_{\fm})\neq
0$, and so by Theorem 2.5, there exists a prime ideal $Q\in
\Assh_{R_{\fm}}M_{\fm}$ such that $\cd_{R_{\fm}}(\fa
R_{\fm},R_{\fm}/Q)=d$. But, then there is exists a prime ideal
$\fp\subseteq \fm$ of $R$ such that $Q=\fp R_{\fm}$. As we have
seen in the proof of Theorem 2.5, $Q\in \Assh_{R_{\fm}}M_{\fm}$,
implies that $$\fp\in \Assh_RM=\Assh_R(R/\fb).$$

Conversely, let $\fm$ be a maximal ideal of $R$ such that there
exists a prime ideal $\fp\in \Assh_R(R/\fb)$ such that
$\fp\subseteq \fm$ and $\cd_{R_{\fm}}(\fa R_{\fm},R_{\fm}/\fp
R_{\fm})=d$. Since $\Var(\fp R_{\fm})\subseteq
\Supp_{R_{\fm}}M_{\fm}$, by [{\bf 9}, Theorem 2.2], it turns out
that $\cd_{R_{\fm}}(\fa R_{\fm},M_{\fm})=d$. But, this implies
that $\fm\in \Supp_RH^d_{\fa}(M)$.

iii) Let $\fm\in \mathcal{P}_{\fa,R/\fb}$. Then by part (ii), we
deduce that $H^d_{\fa R_{\fm}}(M_{\fm})\neq 0$. Hence,  [{\bf 2},
Lemma 2.1] yields that the $R_{\fm}$-module $(H^d_{\fa}(M))_{\fm}$
is not finitely generated. $\Box$

\begin{remark} i) Let $M$ and $N$ be two finitely generated
$R$-modules of finite dimension $d$ so that $\Assh_RN=\Assh_RM$.
Having Theorem 2.5 in mind, it becomes clear that
$\Att_RH^d_{\fa}(N)=\Att_RH^d_{\fa}(M)$. Also, it follows by
Corollary 2.7 (ii) that $\Supp_RH^d_{\fa}(N)=\Supp_RH^d_{\fa}(M)$.
In particular, $H^d_{\fa}(N)=0$ if and only if $H^d_{\fa}(M)=0$.\\
ii) Let $R$ be a ring of finite dimension $d$ and $\fa$ an ideal
of $R$. Also, let $M$ be a finitely generated $R$-module. If $M$
is faithful, then it follows by [{\bf 2}, Theorem 3.3 (b)] that
$\Supp_RH^d_{\fa}(M)=\mathcal{P}_{\fa}$. It is perhaps worth
pointing out that by part (i), this conclusion for $M$ remains
valid under the weaker assumption that $\Assh_RM=\Assh_RR$.
\end{remark}

The following lemma will be needed in the proof of our last result in this section.

\begin{lemma} Let $\fa$ and $\fb$ be two ideals of $R$ and $c$
a natural number. Assume that $M$ is a finitely generated
$R$-module so that $\cd_R(\fa,M)\leq c$. Then there is a natural
isomorphism $$H^c_{\fa}(M/\fb M)\simeq H^c_{\fa}(M)/\fb
H^c_{\fa}(M).$$
\end{lemma}

{\bf Proof.} Let $U=R/\Ann_RM$. Since $\Supp_RU=\Supp_RM$, it
follows by [{\bf 9}, Theorem 2.2], that $H^i_{\fa U}(U)=0$ for all
$i>c$. Hence  $H^i_{\fa U}(\cdot)$ is a right exact functor on the
category of $U$-modules and $U$-homomorphisms. Thus
$$\begin{array}{llll} H^c_{\fa}(M/\fb M) & \simeq
H^c_{\fa U}(U)\otimes_UM/\fb M\\
& \simeq (H^c_{\fa
U}(U)\otimes_UM)\otimes_RR/\fb\\
& \simeq H^c_{\fa}(M)/\fb H^c_{\fa}(M).\Box
\end{array}$$

\begin{theorem} Let $\fa$ be an ideal of $R$ and $M$ a finitely
generated $R$-module such that $c:=\cd_R(\fa,M)\neq -\infty$. Let
$\fW$ be the set of all $\fp\in Supp_RM$ such that $\dim R/\fp=
\cd_R(\fa,R/\fp)=c$ and $\fX:=\fW\cap \Ass_RM$. \\
i) If $\fb:=\cap_{\fp\in \fX}\fp$, then
$\mathcal{P}_{\fa,R/\fb}\subseteq \Supp_RH^c_{\fa}(M)$.\\
ii) If $H^c_{\fa}(M)$ is representable, then $\fX\subseteq
\Att_RH^c_{\fa}(M)$.\\
iii) Assume that $H^c_{\fa}(M)$ is representable. If $\fp\in
\Att_RH^c_{\fa}(M)$ is so that $\dim R/\fp=c$, then $\fp\in \fW$.
\end{theorem}

{\bf Proof.} By [{\bf  1}, p.263, Proposition 4], there is a
submodule $N$ of $M$ such that $\Ass_R(M/N)=\fX$. In particular,
$\dim M/N=c$. Since $\Supp_RN\subseteq \Supp_RM$, by [{\bf  9},
Theorem 2.2], we have $H^i_{\fa}(N)=0$ for all $i>c$. Thus, the
exact sequence $$0\lo N\lo M\lo M/N\lo 0$$ provides the following
exact sequence of local cohomology modules
$$\dots \lo H^c_{\fa}(N)\lo H^c_{\fa}(M)\lo H^c_{\fa}(M/N)\lo 0.$$
Thus $\Supp_RH^c_{\fa}(M/N)\subseteq \Supp_RH^c_{\fa}(M)$, and so
(i) follows by Corollary 2.7 (ii). If $H^c_{\fa}(M)$ is
representable, then the above exact sequence implies that
$\Att_RH^c_{\fa}(M/N)\subseteq \Att_RH^c_{\fa}(M)$, and so (ii)
follows by Theorem 2.5.

Next, we prove (iii). Let $\fp\in \Att_RH^c_{\fa}(M)$ be so that
$\dim R/\fp=c$. By [{\bf 12}, 2.5], there is a submodule $N$ of
$H^c_{\fa}(M)$ such that $\fp=N:_RH^c_{\fa}(M)$. Hence $\fp
H^c_{\fa}(M)\subseteq N$, and so by Lemma 2.9, it turns out that
$H^c_{\fa}(M)/N$ is isomorphic to a quotient of $H^c_{\fa}(M/\fp
M)$. Now, by the Independence Theorem [{\bf 3}, Theorem 4.2.1], we
have the following isomorphisms
$$\begin{array}{llll} H^c_{\fa}(M/\fp M) & \simeq H^c_{\fa R/\fp}
(M/\fp M)\\
& \simeq H^c_{\fa R/\fp}(R/\fp)\otimes_{R/\fp}M/\fp M\\
&\simeq H^c_{\fa}(R/\fp)\otimes_RM. \end{array}$$
Thus
$H^c_{\fa}(M/\fp M)$ is Artinian and $\fp\in \Att_RH^c_{\fa}(M/\fp
M)$. Because, by [{\bf 7}, Corollary 3.3] for an Artinian
$R$-module $A$ and a finitely generated $R$-module $N$, we have
$$\Att_R(A\otimes_RN)=\Att_RA\cap \Supp_RN,$$
the conclusion follows by Theorem 2.5. $\Box$

\section{Lichtenbaum-Hartshorne Vanishing Theorem}

Let the situation be as in Theorem 2.5. In the case that the ideal
$\fa$ is the intersection of finitely many maximal ideals of $R$,
we can find a better description of the set $\Att_RH^d_{\fa}(M)$.
We do this in the next result. The last assertion of this result
might be considered as the generalization of Grothendieck's
non-Vanishing Theorem to semi-local rings.

\begin{proposition} Assume that $\fm_1,\dots ,\fm_t$ are maximal
ideals of $R$ and $M$ a finitely generated $R$-module of finite
dimension $d$. Let $\fa=\bigcap_{i=1}^t\fm_i$. Then
$$\Att_RH^d_{\fa}(M)=\{\fp\in \Assh_RM:\exists 1\leq i\leq t
\ such \ that \ \fp\subseteq \fm_i \ and \
\Ht\frac{\fm_i}{\fp}=d\}.$$ In particular, if $R$ is semi-local
with the only maximal ideals $\fm_1,\dots ,\fm_t$, then
$\Att_RH^d_{\fa}(M)=\Assh_RM$, and so $H^d_{\fa}(M)\neq 0$
whenever $M$ is nonzero.
\end{proposition}

{\bf Proof.} Let $1\leq i\leq t$. Since
$\Supp_RH^d_{\fm_i}(M)\subseteq \{\fm_i\}$, by Lemma 2.1 (ii) and
the Flat Base Change Theorem, it turns out that
$H^d_{\fm_i}(M)\simeq H^d_{\fm_i R_{\fm_i}}(M_{\fm_i})$. Hence,
applying the Mayer-Vietoris sequence for local cohomology provides
the following natural isomorphisms
$$H^d_{\fa}(M)\simeq \displaystyle{\bigoplus_{i=1}^t}H^d_{\fm_i}(M)\simeq
\displaystyle{\bigoplus_{i=1}^t}H^d_{\fm_i
R_{\fm_i}}(M_{\fm_i}).$$ By [{\bf  13}, Theorem 2.2], for a
finitely generated module $M$ over a local ring $(U,\fn)$, we have
$\Att_UH^d_{\fn}(M)=\Assh_UM$, where $d=\dim M$. Thus by Lemma 2.1
(i), we conclude that
$$\begin{array}{llll} \Att_RH^d_{\fa}(M)& =\displaystyle{\bigcup_{i=1}^t}
\{\fp\in \Assh_RM: \fp R_{\fm_i}\in \Assh_{R_{\fm_i}}M_{\fm_i}\
and \ \dim R_{\fm_i}/\fp R_{\fm_i}=d \}\\&=\{\fp\in
\Assh_RM:\exists 1\leq i\leq t \ such \ that \ \fp\subseteq \fm_i
\ and \ \Ht\frac{\fm_i}{\fp}=d\}.\end{array}$$ The last assertion
is immediate by the first one. $\Box$

\begin{remark} Let $A$ be an Artinian $R$-module. Suppose that
$\Supp_RA=\{\fm_1,\dots ,\fm_t\}$ and put $\fM=\cap_{i=1}^t\fm_i$.
Let $T$ denote the $\fM$-adic
completion of $R$.\\
i) Sharp [{\bf  16}] showed that $A$ has a natural structure as a
module over $T$. Let $\theta:R\lo T$ denote the natural ring
homomorphism. The $T$-module structure of $A$ is such that for any
element $r\in R$ the multiplication by $r$ on $A$ has the same
effect as the multiplication of $\theta(r)\in T$. Furthermore a
subset of $A$ is an $R$-submodule of $A$ if and only if it is
a $T$-submodule of $A$.\\
ii) Let $\fa\subseteq \fb$ denote two ideals of $R$ and
$B:=\Sigma_{n\in \mathbb{N}}<\fb>(0:_A\fa^n)$. By [{\bf 10},
Theorem 2.4], the following are equivalent:

a) For any $l\in \mathbb{N}$, there is an integer $n=n(l)$ such
that $0:_A\fM^l\subseteq <\fb>(0:_A\fa^n)$.

b) $B=A$.

c) $\Rad(\fp+\fa T)\varsubsetneq \Rad(\fp+\fb T)$ for all $\fp\in
\Att_TA$.\\
iii) Let $\fa,\fb$ and $B$ be as in (ii), and let $A=S_1+\dots
+S_n$ be a minimal secondary representation of $A$ as a
$T$-module. We can order the elements of $\Att_TA=\{\fp_1,\dots
,\fp_n\}$ such that for an integer $0\leq l\leq n$,
$\Rad(\fp_i+\fa T)\varsubsetneq \Rad(\fp_i+\fb T)$ for all $1\leq
i\leq l$, while $\Rad(\fp_i+\fa T)=\Rad(\fp_i+\fb T)$ for all
$l+1\leq i\leq n$. Then $S_1+\dots +S_l$ is a minimal secondary
representation of $B$. This follows by [{\bf  10}, Theorem 2.8].
Also, it is a routine check to see that
$\Sigma_{i=l+1}^n(S_i+B)/B$ is a minimal secondary representation
of $A/B$ as a $T$-module.
\end{remark}

\begin{theorem} Let $M$ be a finitely generated $R$-module of
finite dimension $d$ and $\mathcal{P}$ a finite subset of $\Max
R$. Let $\fM=\cap_{\fm\in \mathcal{P}}\fm$ and $T$ denote the
$\fM$-adic completion of $R$. \\
i) Let $\fa$ and $\fb$ be two ideals of $R$ such that
$\mathcal{P}_{\fa,M}=\mathcal{P}_{\fb,M}=\mathcal{P}$. If either
$\fa\subseteq \fb$ or $\Att_TH^d_{\fa}(M)\subseteq
\Att_TH^d_{\fb}(M)$, then
$H^d_{\fa}(M)$ is isomorphic to a quotient of $H^d_{\fb}(M)$.\\
ii) Let $\fa$ and $\fb$ be as in (i). If
$\Att_TH^d_{\fa}(M)=\Att_TH^d_{\fb}(M)$, then $H^d_{\fa}(M)\simeq
H^d_{\fb}(M)$. \\
iii) For all ideals $\fc$ of $R$, there are at most
$2^{\mid\Assh_T(M\otimes_RT)\mid}$ non-isomorphic top local
cohomology modules $H^d_{\fc}(M)$ such that
$\Supp_RH^d_{\fc}(M)=\mathcal{P}$.
\end{theorem}

{\bf Proof.} Let $A=H^d_{\fM}(M)$,
$$B_1=\displaystyle{\sum_{n\in \mathbb{N}}}<\fM>(0:_A\fa^n)$$ and
$$B_2=\displaystyle{\sum_{n\in \mathbb{N}}}<\fM>(0:_A\fb^n).$$ Then,
Theorem 2.3 yields the natural isomorphisms $H^d_{\fa}(M)\simeq
A/B_1$ and $H^d_{\fb}(M)\simeq A/B_2$. Let $A=S_1+\dots +S_n$ be a
minimal secondary representation of $A$ as a $T$-module and set
$$\fZ_j :=\Att_TA\setminus \Att_T(A/B_j)$$ for $j=1,2$. Then by
Remark 3.2 (iii), $B_j=\Sigma_{\fp_i\in \fZ_j}S_i$ for $j=1,2$.
Thus, if either $\fa\subseteq \fb$ or $\Att_TH^d_{\fa}(M)\subseteq
\Att_TH^d_{\fb}(M)$, then $B_2\subseteq B_1$, and so
$H^d_{\fa}(M)$ is isomorphic to a quotient of $H^d_{\fb}(M)$.
Also, if $\Att_TH^d_{\fa}(M)=\Att_TH^d_{\fb}(M)$, then $B_1=B_2$,
and so $H^d_{\fa}(M)\simeq H^d_{\fb}(M)$.

Next, we are going to prove (iii). Since $T$ and $\Pi_{\fm\in
\mathcal{P}}\hat{R}_{\fm}$ are isomorphic as $R$-modules, by the
Flat Base Change Theorem, we have the following isomorphisms

$$\begin{array}{llll} H^d_{\fM T}(M\otimes_RT)&
\simeq H^d_{\fM}(M)\otimes_RT\\
& \simeq \displaystyle{\bigoplus_{\fm\in
\mathcal{P}}}(H^d_{\fM}(M)\otimes_R\hat{R}_{\fm})\\
& \simeq \displaystyle{\bigoplus_{\fm\in \mathcal{P}}}(H^d_{\fM
R_{\fm}}(M_{\fm}) \otimes_{R_{\fm}}\hat{R}_{\fm})\\
& \simeq \displaystyle{\bigoplus_{\fm\in
\mathcal{P}}}H^d_{\fm R_{\fm}}(M_{\fm})\\
& \simeq \displaystyle{\bigoplus_{\fm\in \mathcal{P}}}H^d_{\fm}(M)\\
& \simeq H^d_{\fM}(M). \end{array}$$ The last isomorphism follows
by the Mayer-Vietoris sequence for local cohomology. It now is
easy to check that each of these isomorphisms is also a
$T$-isomorphism. Next, as $\fM T$ is the intersection of all
maximal ideals of the semi-local ring $T$, it follows by
Proposition 3.1 that
$$\Att_TH^d_{\fM}(M)=\Att_TH^d_{\fM T}(M\otimes_RT)=
\Assh_T(M\otimes_RT).$$ Now, the claim follows by part (ii).
$\Box$

As an immediate application of Theorem 3.3, we deduce Theorem 1.6
and Proposition 1.5 of [{\bf 5}].

\begin{corollary} Let $\fa$ and $\fb$ be two ideals of a local ring
$(R,\fm)$ and $M$ a finitely generated $R$-module. Let $d=\dim
M$. \\
i) If either $\fa\subseteq \fb$ or
$\Att_{\hat{R}}H^d_{\fa}(M)\subseteq \Att_{\hat{R}}H^d_{\fb}(M)$,
then $H^d_{\fa}(M)$ is isomorphic to a quotient of $H^d_{\fb}(M)$.\\
ii) If $\Att_{\hat{R}}H^d_{\fa}(M)=\Att_{\hat{R}}H^d_{\fb}(M)$,
then $H^d_{\fa}(M)\simeq H^d_{\fb}(M)$.\\
iii) The number of non-isomorphic top local cohomology modules
$H^d_{\fc}(M)$ is at most $2^{\mid \Assh_{\hat {R}}\hat {M}\mid}$
for all ideals $\fc$ of $R$.
\end{corollary}

\begin{example} It might be of some interest to replace
$\hat{R}$ by $R$ in Corollary 3.4 (ii). But, as we show in the
sequel, this would not be the case. To this end, we use an example
of Brodmann and Sharp (see [{\bf 3}, Exercise 8.2.9]). Let $K$ be
a field of characteristic 0. Let $R':=K[X,Y,Z]$, $\fm':=(X,Y,Z)$
and $\fb=(Y^2-X^2-X^3)$. Set $R:=(R'/\fb)_{\fm'/\fb}$ and let
$\fp$ denote the extension of the ideal
$$(X+Y-YZ,(Z-1)^2(X+1)-1)$$ of $R'$ to $R$. As it is mentioned in
[{\bf 3}, Exercise 8.2.9], it follows that $R$ is a 2-dimensional
local domain and that $\fp \hat{R}$ is a prime ideal of $\hat{R}$
with $\dim \hat{R}/\fp \hat{R}=1$. Also, it follows that
$H^2_{\fp}(R)\neq 0$, (see again [{\bf 3}, Exercise 8.2.9]). So
$\Att_{\hat{R}}H^2_{\fp}(R)$ is not empty. Now, let $\fp^*$ be a
minimal associated prime ideal of $\hat{R}$ such that
$\fp^*\subseteq \fp \hat{R}$. Then the inclusion must be strict,
because otherwise we would have
$$\fp=\fp \hat{R}\cap R=\fp^*\cap R\in \Ass_R\hat{R}=\{(0)\},$$ a
contradiction. This yields that $\dim \hat{R}/\fp^*=2$, and so
$\fp^*\in \Assh_{\hat{R}}\hat{R}$. On the other hand, we have
$$\dim \hat{R}/\fp \hat{R}+\fp^*=\dim \hat{R}/\fp \hat{R}=1.$$ Hence
$\fp^*$ does not belong to $\Att_{\hat{R}}H^2_{\fp}(R)$. Thus, if
$\fm$ denotes the maximal ideal of the local ring $R$, then
$$\emptyset \neq \Att_{\hat{R}}H^2_{\fp}(R)\subsetneqq
\Att_{\hat{R}}H^2_{\fm}(R)=\Assh_{\hat{R}}\hat{R}.$$ In
particular, it becomes clear that $H^2_{\fp}(R)$ and
$H^2_{\fm}(R)$ are not isomorphic. On the other hand, we have
$$\Att_{R}H^2_{\fp}(R)=\Att_{R}H^2_{\fm}(R)=\{(0)\}.$$ We therefore
conclude that, it is not possible to replace $\hat{R}$ by $R$ in
Corollary 3.4 (ii).
\end{example}

The following is an analogue of the Lichtenbaum-Hartshorne
Vanishing Theorem for general Noetherian rings.

\begin{theorem} Let $\fa$ be an ideal of $R$ and $M$ a finitely
generated $R$-module of finite dimension $d$. Let
$\fM=\bigcap_{\fm\in \mathcal{P}_{\fa,M}}\fm$ and $T$ denote the
$\fM$-adic completion of $R$. Then the
following  are equivalent:\\
i) $H^d_{\fa}(M)=0$.\\
ii) $H^d_{\fM}(M)=\displaystyle{\sum_{n\in
\mathbb{N}}}<\fM>(0:_{H^d_{\fM}(M)}\fa^n)$.\\
iii) For any integer $l\in \mathbb{N}$, there exists an $n=n(l)\in
\mathbb{N}$ such that $$0:_{H^d_{\fM}(M)}\fa^l\subseteq
<\fM>(0:_{H^d_{\fM}(M)}\fa^n).$$\\
iv) $\dim T/\fa T+\fp>0$ for all
$\fp\in \Assh_T(M\otimes_RT)$.\\
v) $\cd_R(\fa,R/\fp)<d$ for all $\fp\in \Assh_RM$.
\end{theorem}

{\bf Proof.} Let $\fp\in \Assh_T(M\otimes_RT)$. Then, it is easy
to see that $\dim T/\fa T+\fp>0$ if and only if $\Rad(\fp+\fa
T)\varsubsetneq \Rad(\fp+\fM T)$. Therefore, the equivalence of
the conditions (i), (ii) and (iv) is clear by Theorem 2.3 and
Remark 3.2 (ii). Note that in the proof of Theorem 3.3, we have
seen that $\Att_TH^d_{\fM}(M)= \Assh_T(M\otimes_RT)$.

Since $\fa\subseteq \fM$, any element of $H^d_{\fM}(M)$ is
annihilated by some power of $\fa$. Thus $iii)\Rightarrow ii)$
becomes clear.

$ii)\Rightarrow iii)$ Let $A=H^d_{\fM}(M)$ and $l$ a fixed natural
number. Then $(0:_A\fa^l)/<\fM>(0:_A\fa^l)$ is a Noetherian
$R$-module and so the sequence $\{(0:_A\fa^l)\cap
<\fM>(0:_A\fa^n)\}_{n\in \mathbb{N}}$ satisfies the ascending
chain condition. Thus, it follows by [{\bf 10}, Lemma 2.1] that
(ii) implies (iii).

By Grothendieck's Vanishing Theorem, it turns out that
$\cd_R(\fa,R/\fp)\leq d$ for all $\fp\in \Supp_RM$. Therefore, the
equivalence (i) and (v) is immediate by Theorem 2.5. $\Box$


\end{document}